\documentclass{elsart3-1}
 \usepackage{graphicx}
\usepackage{amssymb}
\usepackage[latin1]{inputenc}

\usepackage[english,francais]{babel}
\newtheorem{theo}{Theorem}[section]

\newtheorem{defi}[theo]{Definition\rm}

\newtheorem{example}{\it Example\/}
\newcommand{\D}{\mathcal D}

\renewcommand{\theequation}{\arabic{equation}}

\newcommand{\rp}{\mathbb R\mathbb P}
\newcommand{\cp}{\mathbb C\mathbb P}
\newcommand{\PPP}{\mathcal P}
\newcommand{\vrt}{\operatorname{Vert}}
\newcommand{\edge}{\operatorname{Edge}}

\newcommand{\N}{\mathbb N}

\newcommand{\C}{\mathbb C}
\newcommand{\R}{\mathbb R}
\newcommand{\dv}{\operatorname{div}}

\renewcommand{\setminus}{\smallsetminus}
\def\og{\leavevmode\raise.3ex\hbox{$\scriptscriptstyle\langle\!\langle$~}}
\def\fg{\leavevmode\raise.3ex\hbox{~$\!\scriptscriptstyle\,\rangle\!\rangle$}}

\begin{document}
\begin{frontmatter}


\selectlanguage{english}
\title{Enumeration of curves via floor diagrams}


\selectlanguage{english}
\author[authorlabel1]{Erwan Brugall\'e},
\ead{brugalle@math.jussieu.fr}
\author[authorlabel2]{Grigory Mikhalkin}
\ead{mikha@math.toronto.edu}

\address[authorlabel1]{Universit\'e Pierre et Marie Curie, Institut
  Mathématiques de Jussieu, 175 rue du Chevaleret, 75 013 Paris, France}
\address[authorlabel2]{ University of Toronto, Department of Mathematics,
40 St. George St., Toronto, Ontario M5S 2E4, Canada}


\begin{abstract}
\selectlanguage{english}
In this note we compute some enumerative invariants of real and
complex projective spaces by means of some enriched graphs called
floor diagrams.

\end{abstract}
\end{frontmatter}

\selectlanguage{english}
\section{Introduction: enumerative invariants of real and complex projective
  spaces}\label{intro}
The subject of this note is the number of curves passing through
configuration of linear spaces in $\rp^n$ and $\cp^n$. To set up
an enumerative problem we fix integer numbers $d\ge 1$, $g\ge 0$ and
$n\ge 2$. Then we look at the algebraic curves of degree $d$ and
genus $g$ in $\rp^n$ or $\cp^n$. Namely, by the curves of genus $g$
in $\cp^n$ we mean the images of Riemann surfaces of genus $g$ under
holomorphic maps to projective spaces; by the curves in $\rp^n$ we
mean the real points of those curves in $\cp^n$ that are invariant
under the involution of complex conjugation.

By the Riemann-Roch formula, the space of all such curves has
dimension greater or equal than  $(n+1)d + (n-3)(1-g)$. Furthermore,
we always have equality in the case when $n=2$ or in the case when
$g=0$ for $n>2$. These are the two cases that we are concerned with
in this note.

Let us fix a generic configuration $\PPP$
of projective-linear subspaces
(of different dimensions) in $\rp^n$ or $\cp^n$, so that
we have $l_j$ subspaces of dimension $j=0,\dots,n-2$.
These subspaces are known as {\em constraints of dimension $j$}.
It can be shown that if
\begin{equation}\label{pl}
\sum\limits_{j=0}^{n-2}l_j(n-1-j)=(n+1)d+(n-3)(1-g)
\end{equation}
then
the number of curves of degree $d$ and genus $g$ passing via $\PPP$
is finite (recall that our assumption was that $g=0$ whenever $n>2$).
Furthermore, if we work in $\cp^n$ then this number does not depend
on the choice of $\PPP$ and is known as the Gromov-Witten number,
which we denote by $N_{d,g}^{(n)}(l_0,\dots,l_{n-2})$.
These invariants are well-known; they
were computed by Kontsevich (see \cite{KonMan1}) in the case of $g=0$
and arbitrary $n$ and by Caporaso and Harris (see \cite{CapHar1}) in the case of $n=2$
and arbitrary $g$.

Much less is known if we work in $\rp^n$.
Certainly in this case the number of curves
of degree $d$ and genus $g$ depends not only
on the numbers $l_j$ but also on the
choice of the configuration $\PPP$. Nevertheless,
Welschinger proved in \cite{Wel1}, \cite{Wel2} that in the case when  $n=2,3$, $g=0$ and $l_j=0$
for any $j>0$,
 there is a consistent choice of signs $\pm 1$ such that the number
of the corresponding real curves counted with the sign depends only
on $d$ (note that the assumption $l_j=0$ for $j>0$
allows us to determine $l_0$ once we fix $d$ by $(\ref{pl})$ ).
These numbers are known as the Welschinger
numbers, we denote them with $W_d^{(n)}$.

The technique for computation of the numbers $W_d^{(2)}$ was
given in \cite{Mik1} as an application of tropical geometry.
It was used in \cite{IKS1} to produce several estimates for these numbers.
In particular, we have $W_d^{(2)}>0$ for any $d$.
By symmetry reason we have $W_d^{(3)}=0$ for any even $d$,
but the values of $W_d^{(3)}$ remained unknown for odd $d$ even for $d=5$.

In this note we announce a technique which allows a simultaneous
computation of $N_{d,g}^{(n)}(l_0,\dots,l_{n-2})$
and $W_d^{(n)}$ (for $n=2,3$).
The details will be given in \cite{BM} along with a somewhat
more general computation of the number of real curves in the case
when the configuration $\PPP$ contains pairs of complex conjugate subspaces.

\section{Floor diagrams}
Let $\Gamma$ be a finite oriented graph. We say that $\Gamma$ is
{\em acyclic} if it does not contain any non-trivial {\em oriented}
cycle (in the same time the first Betti number of $\Gamma$ may be
high). We denote the set of its vertices with
$\overline{\vrt}(\Gamma)$ and the set of its (open) edges with
$\edge(\Gamma)$. We denote by $\vrt^\infty(\Gamma)$ the set of
sinks (i.e. the vertices such that all their adjacent edges are incoming),
and  with $\edge^\infty(\Gamma)$ the set of edges adjacent
to a sink. Finally we put
$\vrt(\Gamma)=\overline{\vrt(\Gamma)}\setminus\vrt^\infty(\Gamma)$.

We say that $\Gamma$ is a {\em weighted graph} if each
edge of $\Gamma$ is prescribed a natural weight, i.e. we are given
a function $w:\edge(\Gamma)\to\N$.
The weight allows one to define the {\em divergence} at the vertices.
Namely, for a vertex $v\in\vrt(\Gamma)$
we define the divergence $\dv(v)$ to be the sum of the weights of all
outgoing edges minus the sum of the weights of all incoming edges.


\begin{defi}
A connected weighted oriented graph $\D$ is called a
{\em floor diagram of genus $g$ and degree $d$}
if the following conditions hold.
\begin{itemize}
\item The oriented graph $\D$ is acyclic.
\item We have $\dv(v)>0$ for any $v\in\vrt(\D)$ and $\dv(v)=-1$
for every $v\in \vrt^{\infty}(\D)$.
\item The first Betti number $b_1(\D)$ equals $g$.
\item 
The set $\vrt^{\infty}(\D)$ consists of $d$ elements.
\end{itemize}
\end{defi}
We call a vertex $v\in\vrt(\D)$ {\em a floor of degree $\dv(v)$}.

Let $l_0,\dots,l_{n-2}\ge0$ be integer numbers subject to $(\ref{pl})$
and $\PPP=\{x^{(0)}_1,\dots,x^{(0)}_{l_0},\dots
x^{(n-2)}_1,\dots,x^{(n-2)}_{l_{n-2}}\}$ be the ordered set
of $\sum\limits_{j=0}^{n-2}l_j$ elements.
We define $\dim(x^{(j)}_k)=j$.
Let $m:\PPP\to\D\setminus\vrt^{\infty}(\D)$ be a map.
Let $\D_{>e}$ be the component of $\D\setminus e$ that contains the
arrowhead of $e$.
We define the {\em height} of $e$ by $h(e)=0$ if $g>0$ and otherwise by
$$h(e)=\sum_{q\in\PPP,\ m(q)\in\D_{>e}}\big(n-1-\dim(q)\big)
 +1 - w(e)-(n+1)\sum_{v\in\vrt(\D)\cap\D_{>e}}\dv(v).$$
(The meaning of the height is the dimension of the constraint we have to put on $e$
to make $\D_{>e}$ rigid.)
\begin{defi}
The map $m$
is called the {\em marking} of $\D$ if it satisfies to the following
properties.
\begin{itemize}
\item If $m(q)=m(q')$ for $q< q'\in\PPP$ then $m(q)$ is a vertex and
  $\dim(q)>0$.
\item For any $v\in \vrt(\D)$, there exists $q\in\PPP$ such that
  $m(q)=v$.
\item  If $q'>q$ and  $m(q')<m(q)$, then
there exists $q''$ such that $m(q)=m(q'')$ and
$q''>q'$.
\end{itemize}
\end{defi}
A floor diagram $\D$ enhanced with a marking is called the {\em
  marked
floor diagram}, and is said to be {\em marked by $\PPP$}. Two marked
floor diagrams are called {\em equivalent} if they can be identified
by a homeomorphism of oriented graphs. They are called to be of the
{\em same combinatorial type} if there exists a bijection
$\sigma:\PPP\to\PPP$ that preserves the dimension of the constraints
and such that $(\D,m)$ is equivalent to $(\D,m\circ\sigma)$.

To each marked floor diagram $\D$, we may associate
its complex and real {\em multiplicities} $\mu^{\C}(\D)$ and $\mu^{\R}(\D)$.
These multiplicities
record the number of complex and real curves (respectively)
encoded by the diagram.

The definition is inductive by $n$. Take a floor $v\in\vrt(\D)$ of
degree $\dv(v)$.
For the maximum element $x_k^{(j)}$ of $m^{-1}(v)$ (recall that $\PPP$
is ordered), take a linear space
of dimension $j$.
For each other element $x_{k'}^{(j')}$ in $m^{-1}(v)$, take a linear space
of dimension $j'-1$.
Take a linear space of dimension $h(e)$ for each edge $e$ incoming to $v$.
For each edge $e$ outgoing from $v$ take a linear space of dimension
$$n-1-h(e)-\sum_{q\in\PPP,\ m(q)\in e}\big(n-1-\dim(q)\big).$$
If any of these numbers is outside of the range between
0 and $(n-2)$ then we set both the real and complex multiplicity of
$\D$ equal to 0.
Otherwise denote the number of resulting $j$-dimensional linear spaces
with $l^{(v)}_j$
and define
$\mu^{\C}(v)=\dv(v)^{l^{(v)}_{n-2}}N^{(n-1)}_{\dv(v),0}(l_0^{(v)},\dots,l^{(v)}_{(n-3)})$
and $$\mu^{\C}(\D,m)=\prod\limits_{v\in\vrt(\D)}\mu^{\C}(v)
\prod\limits_{e\in\edge(\D)}(w(e))^{1+\#(m^{-1}(e))}.$$
Here $\#(m^{-1}(e))$ is the number of elements in the inverse image
of an {\em open} edge $e$.
To fix the base of induction for $n=1$ we define $N^{(1)}_{1,0}=1$ and
$N^{(1)}_{a,b}=0$ for $(a,b)\neq (1,0)$.
If $g=0$ and $l_j=0$ for $j\ge 1$ then we also define the real multiplicity
$$\mu^\mathbb R(\D)=\prod_{v\in\vrt(\D)} W^{(n-1)}(\dv(v))\ \  if \ \
w(\edge(\D))\cap 2\mathbb N=\emptyset \ \ \ and \ \ \ \mu^\mathbb
R(\D)=0 \ otherwise. $$
We put $W^{(1)}_1=1$ and $W^{(1)}_d=0$ for $d>1$.

Note that for $n=2$, both real and complex
multiplicities have a very simple form. For any diagram with
the non-zero multiplicity we have $h(e)=0$ for any
$e\in\edge(\D)$ and $\dv(v)=1$ for any $v\in\vrt(\D)$, so $\mu^\mathbb C(\D)$
is the square of the product of the weight of all edges of $\D$, and
$\mu^\mathbb R(\D)=\mu^\mathbb C(\D)\ mod \ 2$.
Note also that all marked floor diagrams of the same combinatorial type have
the same multiplicities.

In the examples of combinatorial types of marked floor
diagrams we give below, the convention we use  are the following~: the
graph $\D$ is
oriented from up to down, vertices
are represented by ellipses with its degree written inside,
edges are represented by solid lines close of which is written its
weight (except if it is 1),
and the images of $m$ are represented by a set of points on
$\D$. Below any combinatorial type, we write the number of marked
floor diagram of this combinatorial type, and the complex and real (if
any) multiplicities of such a floor diagram.

\begin{example}
In Figure \ref{n=2,d=3}, we depict all $2$-dimensional marked floor
diagrams of degree 3 with non null multiplicity.
In Figure \ref{n=3,d=5} (resp.  \ref{n=3,d=2}), we depict all $3$-dimensional marked floor diagrams of degree 5
  and marked by 10 points (resp. of degree 2
  and marked by 8 lines) with non null multiplicity.

\begin{figure}[h]
\centering
\begin{tabular}{ccccccc}
\includegraphics[width=0.8cm, angle=0]{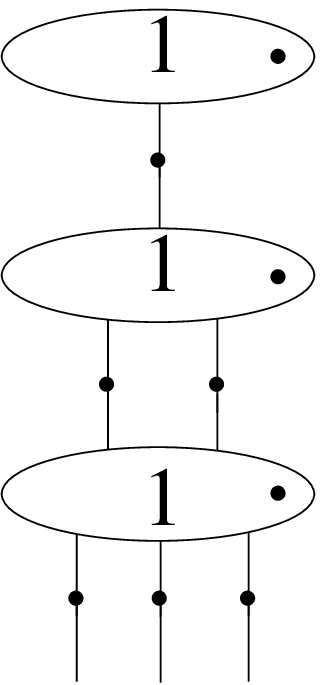}&
 \hspace{4ex}  &
\includegraphics[width=0.8cm, angle=0]{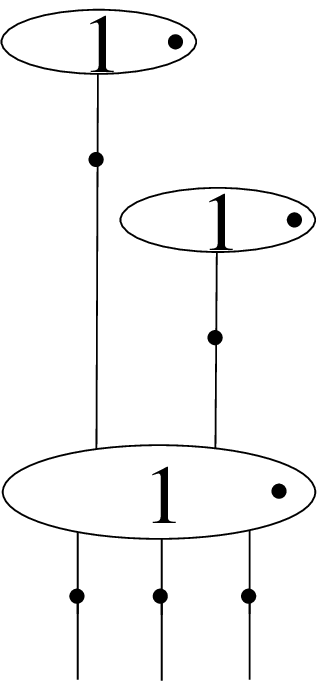}&
 \hspace{4ex}  &
\includegraphics[width=0.8cm, angle=0]{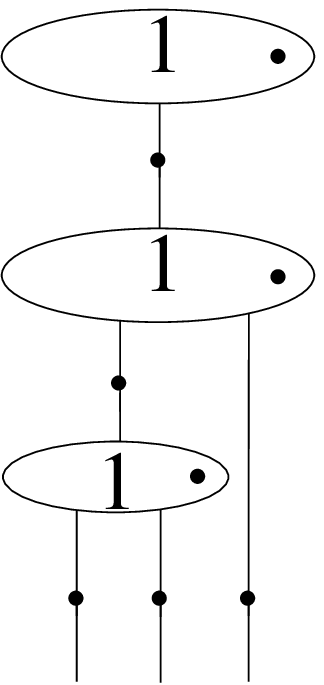}&
 \hspace{4ex}  &
\includegraphics[width=0.8cm, angle=0]{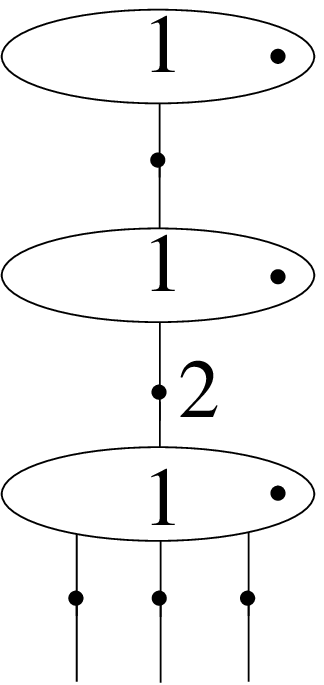}
\\  1, $\mu^\mathbb C=1$ &&  3, $\mu^\mathbb C=1$,
$\mu^\mathbb R=1$  &&  5, $\mu^\mathbb C=1$, $\mu^\mathbb
R=1$
&& 1,  $\mu^\mathbb C=4$, $\mu^\mathbb R=0$
\end{tabular}
\caption{$2$-dimensional marked floor diagrams of degree 3}
\label{n=2,d=3}
\end{figure}
\end{example}

\begin{figure}[h]
\centering
\begin{tabular}{ccccc}
\includegraphics[width=0.9cm, angle=0]{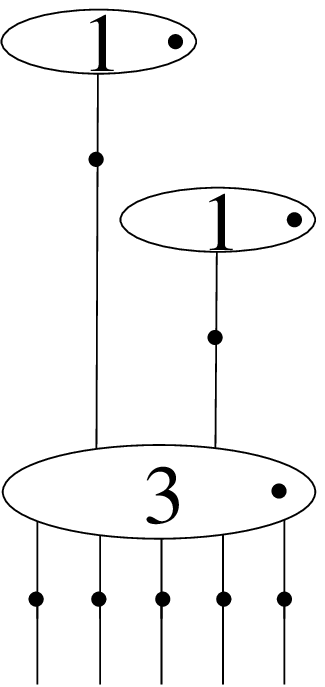}&
 \hspace{4ex}  &
\includegraphics[width=0.9cm, angle=0]{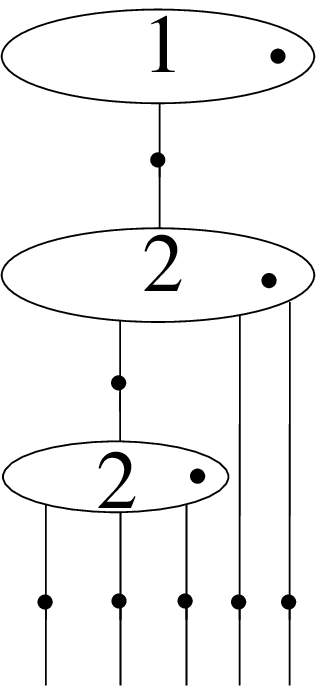}&
 \hspace{4ex}  &
\includegraphics[width=0.9cm, angle=0]{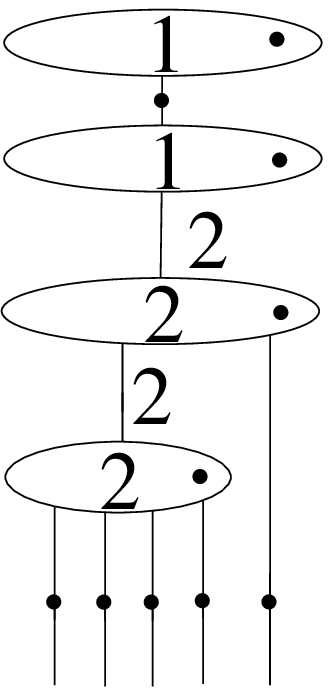}
\\  3, $\mu^\mathbb C=12$, $\mu^\mathbb R=8$ &&  21,
$\mu^\mathbb C=1$,
$\mu^\mathbb R=1$  &&   6, $\mu^\mathbb C=8$, $\mu^\mathbb
R=0$
\end{tabular}
\caption{$3$-dimensional marked floor diagrams of degree 5, genus 0,
  and marked by 10 points}
\label{n=3,d=5}
\end{figure}

\begin{figure}[h]
\centering
\begin{tabular}{ccccccccccccc}
\includegraphics[width=1.0cm, angle=0]{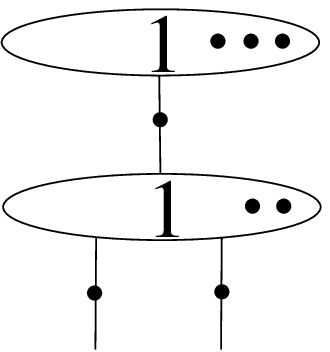}&
 \hspace{2ex}  &
\includegraphics[width=1.0cm, angle=0]{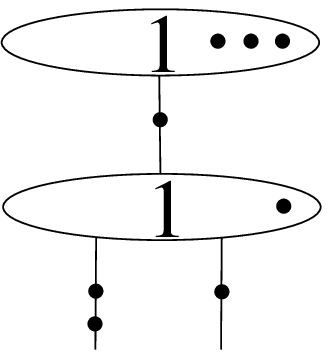}&
 \hspace{2ex}  &
\includegraphics[width=1.0cm, angle=0]{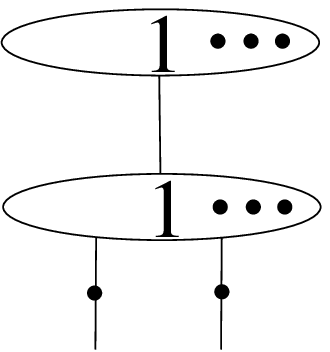}&
 \hspace{2ex}  &
\includegraphics[width=1.0cm, angle=0]{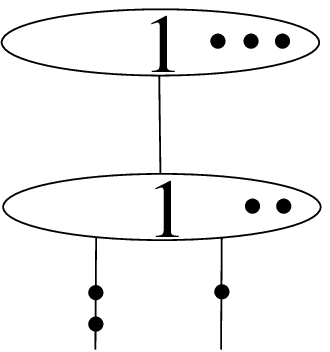}&
 \hspace{2ex}  &
\includegraphics[width=1.0cm, angle=0]{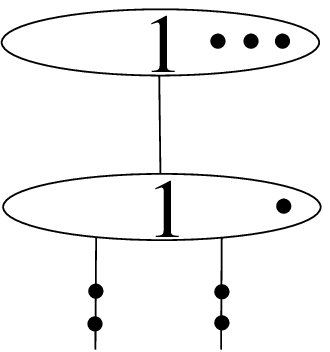}&
 \hspace{2ex}  &
\includegraphics[width=1.0cm, angle=0]{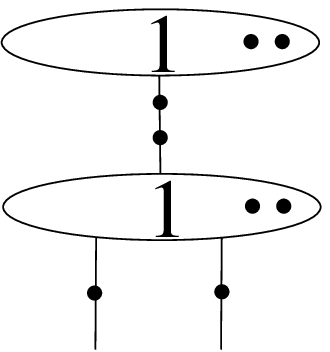}&
 \hspace{2ex}  &
\includegraphics[width=1.0cm, angle=0]{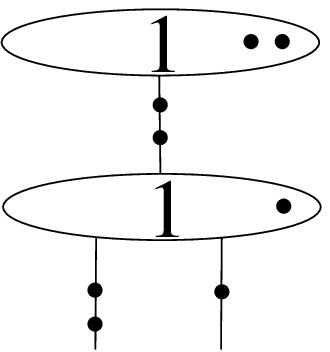}
\\  5, $\mu^\mathbb C=1$ &&   3, $\mu^\mathbb C=1$
  &&   10, $\mu^\mathbb C=1$
&&12,  $\mu^\mathbb C=1$
&&3,  $\mu^\mathbb C=1$
&&5,  $\mu^\mathbb C=1$
&&3,  $\mu^\mathbb C=1$
\end{tabular}

\vspace{2ex}

\begin{tabular}{ccccccccccc}
\includegraphics[width=1.0cm, angle=0]{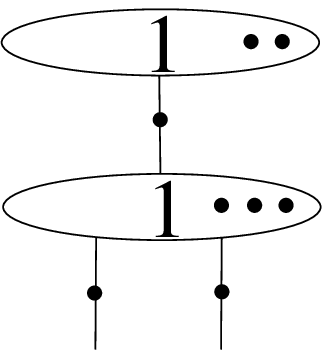}&
 \hspace{2ex}  &
\includegraphics[width=1.0cm, angle=0]{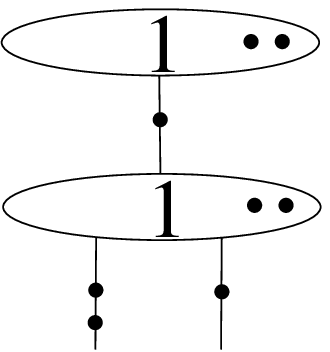}&
 \hspace{2ex}  &
\includegraphics[width=1.0cm, angle=0]{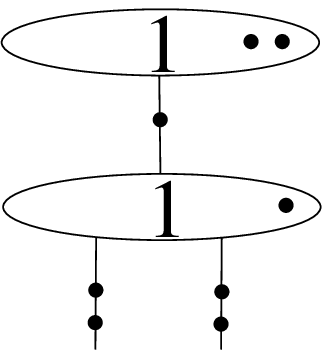}&
 \hspace{2ex}  &
\includegraphics[width=1.0cm, angle=0]{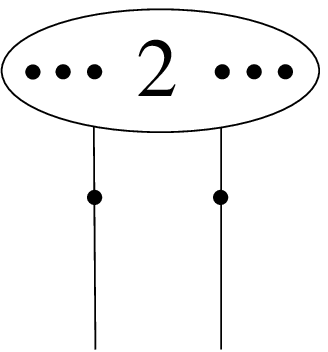}&
 \hspace{2ex}  &
\includegraphics[width=1.0cm, angle=0]{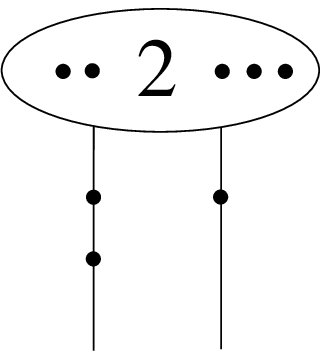}&
 \hspace{2ex}  &
\includegraphics[width=1.0cm, angle=0]{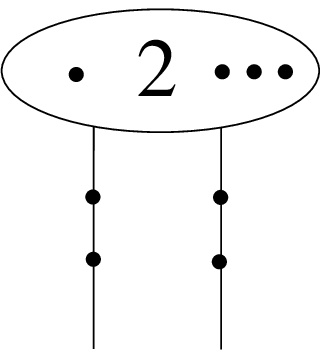}
\\  10, $\mu^\mathbb C=1$ &&   12, $\mu^\mathbb C=1$
  &&   3, $\mu^\mathbb C=1$
&&1,  $\mu^\mathbb C=8$
&&3,  $\mu^\mathbb C=4$
&&3,  $\mu^\mathbb C=2$
\end{tabular}

\caption{$3$-dimensional marked floor diagrams of degree 2, genus 0,
  and marked by 8 lines}
\label{n=3,d=2}
\end{figure}

\section{Main Formulas}\label{results}


\begin{thm}\label{GWform}
For $n=2$ or $g=0$, the number $N^{(n)}_{d,g}(l_0,\dots,l_{n-2})$ is equal
to the sum of the complex multiplicity of all floor
diagrams  of degree $d$, genus $g$ marked by $\PPP$.
\end{thm}
\begin{example}
Using marked floor diagrams depicted in Figures \ref{n=2,d=3},
\ref{n=3,d=5}, and  \ref{n=3,d=2}, we verify that $N^{(2)}_{3,1}(9)=1$,
$N^{(2)}_{3,0}(8)=12$,  $N^{(3)}_{5,0}(10,0)=105$, and   $N^{(3)}_{2,0}(0,8)=92$.
\end{example}
\begin{example}
Using Theorem \ref{GWform}, one can compute the numbers
$N^{(2)}_{d,\frac{(d-1)(d-2)}{2}-1}(\frac{d(d+3)}{2}-1)=3(d-1)^2$ (the degree of the
discriminant) and
$N^{(2)}_{d,\frac{(d-1)(d-2)}{2}-2}(\frac{d(d+3)}{2}-2)=\frac{3}{2}(d-1)(d-2)(3d^2 -
3d-11)$, cf. \cite{Vein}.
\end{example}


\begin{thm}\label{wel}
For $n=2$ or $n=3$, the Welschinger invariant $W^{(n)}_d$ is equal
to $(-1)^{\frac{n(d-1)(d-2)}{2}}$ times the sum of the real
multiplicity of all  $n$-dimensional floor diagrams  of degree $d$,
genus $0$ marked by $\PPP$ with $l_j=0$ for $j>0$.
\end{thm}

The two-dimensional Welschinger invariants were studied by Itenberg,
 Kharlamov and  Shustin \cite{IKS1} who obtained a number of estimates on these
 invariants.
Theorem \ref{GWform} and \ref{wel} allow one not only to recover these
results but also to show that similar results hold for the 3-dimensional
invariants $W^{(3)}$.

\begin{example}
We list the values of
$N^{(3)}_{d,0}(2d,0)$ and $W^{(3)}_d$ up to degree 7.
\begin{center}
\begin{tabular}{ c|c|c|c|c|c|c|c}
 $d$      &\ \  \ \ $1$\ \ \ \  &\ \  \ \  $2$\ \  \ \
 &\ \  \ \ $3$ \ \  \ \ &\ \  \ \  $4$\ \  \ \ &\ \  \ \  $5$\ \  \ \
 &\ \  \ \  $6$ \ \  \ \  &\ \  \ \  $7$\ \  \ \  \\
\hline
 $N^{(3)}_{d,0}(2d,0)$ & $1$ &$0$  &$1$ &$4$ &$105$&$2576$&$122129$
\\ \hline
$W^{(3)}_d$ & $1$ & $0$ & $-1$ & $0$ & $45$ & $0$ & $-14589$
\end{tabular}
\end{center}
\end{example}

\begin{prop}\label{increase n3}
The 3-dimensional Welschinger invariants have the following properties :
\begin{itemize}
\item For any $k> 1$, $|W^{(3)}_{2k+1}|> |W^{(3)}_{2k-1}|$.
\item The sequences $(W^{(3)}_{(2k+1)})$ and $(N^{(3)}_{2k+1,0}(4k+2,0))$ are logarithmically equivalent,
i.e.
$$ \log(|W^{(3)}(2k+1)|)\sim 4k \log k\sim  \log(N^{(3)}_{2k+1,0}(4k+2,0)).$$
\item For any $d\ge 1$, $|W^{(3)}_d|=N^{(3)}_{d,0}(2d,0) \ \ \ mod \ 4$.

\end{itemize}
\end{prop}

To deduce Theorems \ref{GWform} and \ref{wel} we use tropical
geometry. The main technical ingredient is the following lemma
(stated in the notation of Section \ref{intro}, and where $\PPP$ is
here a
configuration of tropical linear subspaces in $\R^n$).

\begin{lem}
 Let $HC$ be a hypercube containing
all vertices of all elements of $\PPP$. Then $HC$ contains all the
vertices of any tropical curve of degree $d$ and genus $g$ in $\mathbb
R^n$ passing
through elements of $\PPP$.
\end{lem}

\section{Further remarks}
\begin{itemize}
\item Results of this note can be extended from projective spaces to
some other (but not all) toric varieties.
\item As in \cite{Mik1}, one can count how many real curves pass
  through some special real configurations $\PPP$. For example, one
  can choose
  8 real lines in $\rp^3$ such that the 92 conics passing
  through these lines are real.
\end{itemize}




We are grateful to Michel Coste, Sergey Fomin, Ilia Itenberg, Viatcheslav
Kharlamov and Oleg Viro for
fruitful discussions. We would like to acknowledge the support from PREA.
The second author is also supported in part by NSERC and the Canada Research Chair
program.


\begin{thebibliography}{00}



\bibitem{BM} Brugall\'e, E. and Mikhalkin, G., Floor decomposition of
  tropical curves, in preparation.
\bibitem{CapHar1}Caporaso, L. and Harris, J., Counting plane curves of
   any genus, Invent. Math., 1998, vol 131, p 345-392.

\bibitem{IKS1}Itenberg, I. and  Kharlamov, V. and Shustin, E.,
  Welschinger invariant and enumeration of real rational curves,
  Int. Math. Research Notices, 2003, vol. 49, p 2639-2653.




\bibitem{KonMan1} Kontsevich, M. and Manin, Yu., Gromov-Witten
              classes, quantum cohomology, and enumerative
              geometry, Comm. Math. Phys., 1994, vol 164, p 525-562.

\bibitem{Mik1}Mikhalkin, G., Enumerative tropical algebraic geometry
  in $\mathbb R^2$, J. Amer. Math. Soc., 2005, vol. 18, p 313-377.

\bibitem{Vein} Vainsencher, I., Enumeration of $n$-fold tangent
  hyperplanes to a surface,
  J. Algebraic Geom.,  1995,  vol 4, p 503-526.

\bibitem{Wel1} Welschinger, J. Y., Invariants of real symplectic
  4-manifolds and lower bounds in real enumerative geometry,
  Invent. Math., 2005, vol 162, p 195-234.

\bibitem{Wel2} Welschinger, J. Y., Spinor states of real rational
              curves in real algebraic convex  3-manifolds and
              enumerative invariants, Duke Math. J., 2005, vol 127, p 89-121.


\end{thebibliography}
\end{document}